\documentclass[letterpaper, 10 pt, conference]{ieeeconf}  

\IEEEoverridecommandlockouts                              
\usepackage{amsmath} 
\usepackage{amsfonts,bm}
\usepackage{mathtools}
\usepackage{subcaption}
\usepackage{algorithm} 
\usepackage{siunitx}
\usepackage[noend]{algpseudocode}

\usepackage{tikz}
\usepackage{tikz-imagelabels}
\usepackage{pgfplots} 
\usepackage{pgfgantt}
\usepackage{pdflscape}
\pgfplotsset{compat=newest} 
\pgfplotsset{plot coordinates/math parser=false} 
\newlength\fwidth
\newlength\fheight
\newtheorem{example}{Example}

\usepackage[colorlinks,citecolor=green,urlcolor=blue,bookmarks=false,hypertexnames=true]{hyperref} 
\newcommand{\dt}{\Delta t}

\newcommand{\q}{\mathbf q}
\newcommand{\p}{\mathbf p}

\newcommand{\vb}{\mathbf b}
\newcommand{\vc}{\mathbf c}
\newcommand{\mA}{\mathbf A}

\newcommand{\N}{\mathcal{N}}
\newcommand{\param}{\pmb \theta}

\newcommand{\reals}{\mathbb{R}}

\newcommand{\thetsig}{\param_{\Sigma}}
\newcommand{\thetgam}{\param_{\Gamma}}
\newcommand{\thetdyn}{\param_{\Psi}} 
\newcommand{\thetobs}{\param_{h}}
\newcommand{\like}{\mathcal{L}}
\newcommand{\probd}{\pi}

\newcommand{\numSteps}{T}

\newcommand{\dimx}{{d}}
\newcommand{\dimpar}{{d_{\param}}}
\newcommand{\dimy}{{d_y}}
\newcommand{\numObs}{n}
\newcommand{\numTraj}{M}
\newcommand{\numBasis}{N}
\newcommand{\numMCMC}{P}
\newcommand{\yn}{\mathcal{Y}_{\numObs}}
\newcommand{\ynm}{\mathcal{Y}_{\numObs}^{(m)}}
\newcommand{\yk}{\mathcal{Y}_{k}}
\newcommand{\X}{\bm{X}}
\newcommand{\Y}{\bm{Y}}
\newcommand{\xn}{\mathcal{X}_{\numObs}}
\usepackage[textsize=tiny]{todonotes}

\newcommand{\Ra}[1]{\color{black} {#1}}
\newcommand{\Rb}[1]{\color{black} {#1}}
\title{\LARGE \bf
Bayesian Identification of Nonseparable Hamiltonian Systems \\ Using Stochastic Dynamic Models
}

\author{Harsh Sharma$^{1*}$, Nicholas Galioto$^{2*}$, Alex A. Gorodetsky$^{2}$, and Boris Kramer$^{1}$
\thanks{$^*$These authors contributed equally.}
\thanks{$^{1}$H. Sharma and B. Kramer are with the Department of Mechanical and Aerospace Engineering, University of California San Diego, San Diego, California, USA. {\tt \small \{hasharma, bmkramer\}@ucsd.edu} }
\thanks{$^{2}$N. Galioto and A.A. Gorodetsky are with the Department of Aerospace Engineering, University of Michigan, Ann Arbor, Michigan, USA.  {\tt \small \{ngalioto, goroda\}@umich.edu}}
}

\begin{document}

\maketitle
\thispagestyle{empty}
\pagestyle{empty}

\begin{abstract}

This paper proposes a probabilistic Bayesian formulation for system identification (ID) and estimation of nonseparable Hamiltonian systems using stochastic dynamic models.
Nonseparable Hamiltonian systems arise in models from diverse science and engineering applications such as astrophysics, robotics, vortex dynamics, charged particle dynamics, and quantum mechanics.
The numerical experiments demonstrate that the proposed method recovers dynamical systems with higher accuracy and reduced predictive uncertainty compared to state-of-the-art approaches. 
The results further show that accurate predictions far outside the training time interval in the presence of sparse and noisy measurements are possible, which lends robustness and generalizability to the proposed approach. A quantitative benefit is prediction accuracy with less than 10\% relative error for more than 12 times longer than a comparable least-squares-based method on a benchmark problem.

\end{abstract}

\section{Introduction}
{\Rb{Nonseparable Hamiltonian systems arise as models in many science and engineering applications such as multibody dynamics and control in robotics~\cite{serra2019control}, the Kozai-Lidov mechanism in astrophysics~\cite{li2014chaos}, particle accelerators in accelerator physics~\cite{forest2006geometric}, 3D vortex dynamics in fluid mechanics~\cite{salmon1988hamiltonian}, and the nonlinear Schr\"odinger equation in quantum mechanics~\cite{colliander2010transfer}.}} These systems demonstrate complex nonlinear behavior while possessing an underlying highly structured geometry encoded by a Hamiltonian. Uncovering a system’s Hamiltonian can reveal key insights into its physical properties such as mass or energy conservation. Learning Hamiltonian models directly from data is becoming increasingly important in diverse areas such as astrophysics, robotics, fluid dynamics, plasma physics, and quantum mechanics where first-principle modeling can yield highly complex model structures or such models are not available.

For these purposes, many recent methods have embedded strong physics-motivated inductive priors into their learning framework to develop structure-preserving neural networks for Hamiltonian systems, e.g., \cite{greydanus2019hamiltonian,chen2019symplectic,zhu2020deep,jin2020sympnets}. In another research direction, techniques based on sparse identification of nonlinear dynamics (SINDy) \cite{chu2020discovering} and  orthogonal polynomials \cite{wu2020structure} have also been used for learning Hamiltonian systems from data, but these techniques tend to break down when the data are noisy/sparse since they rely on numerical approximations of the time derivatives of the data. To handle uncertainty and increase robustness, Bayesian inference techniques based on Gaussian process regression~\cite{bertalan2019learning,offen2022symplectic} and/or consideration of stochastic dynamics~\cite{galioto2020hamiltonian} have been developed. However, a majority of these approaches assume that the Hamiltonian system is separable, i.e., the system Hamiltonian can be written as $H(\q,\p)=T(\p) + U(\q)$ where $\q$ is the position, $\p$ is the momentum,  $T(\p)$ is the kinetic energy, and $U(\q)$ is the potential energy.
    
For nonseparable Hamiltonian systems, nonseparable symplectic neural networks were recently developed in \cite{xiong2020nonseparable} by embedding an appropriate symplectic integrator into the neural network architecture. Instead of learning the symplectic maps, generating function neural networks (GFNNs) in \cite{chen2021data} learn Hamiltonian models by approximating the generating functions corresponding to these symplectic maps. Although both of these approaches learn nonseparable Hamiltonian systems from noisy data, they use an optimization objective and modeling format that requires large datasets. These methods use data from over a thousand (up to a million in the case of GFNN) short trajectories, each consisting of five or fewer data points, even for one-dimensional problems.

In this paper, we learn a nonseparable Hamiltonian system from noisy and sparse data using a stochastic model to account for model errors that always exist in any approximation format. Previously, we have shown the consideration of such stochastic process noise aids system recovery by smoothing the learning objective~\cite{galioto2020bayesian}. Building on the Bayesian {\Rb{system identification (ID)}} work in \cite{galioto2020bayesian}, we embed the physics underlying nonseparable Hamiltonian structure into the learning formulation to develop a probabilistic learning method that preserves the symplectic structure intrinsic to Hamiltonian dynamics. The main contributions of this work are:
    \begin{enumerate}
        \item We extend the Bayesian system ID method of \cite{galioto2020bayesian} to nonseparable Hamiltonian systems by embedding an explicit structure-preserving numerical integrator within the physics-informed objective;
        \item We present detailed numerical results for direct comparison between symplectic and non-symplectic approaches that demonstrate the advantage of symplectic structure preservation. The symplectic approach achieves $30\%$ reduction in relative state error over the non-symplectic approach;
        \item We apply the Bayesian method to a dataset with multiple trajectories and non-Gaussian measurement noise for the first time to demonstrate its versatility and robustness. {\Rb{Through this study, we also show that the Gaussian filtering approaches we leverage are not overly restrictive even for nonlinear dynamics}}. After training on a noisy dataset, the proposed method has a relative error below 10\% and does so for more than 12 times longer than the least-squares comparison method on a trajectory outside the training set.
    \end{enumerate}

This paper is structured as follows. Section \ref{sec: background} reviews the basics of nonseparable Hamiltonian systems and describes a probabilistic formulation of the system ID problem. Section \ref{sec: learning} presents the proposed structure-preserving algorithm for learning nonseparable Hamiltonian systems. In Section \ref{sec: numerical}, we apply the proposed Bayesian algorithm to two datasets generated from nonseparable Hamiltonian systems. We also compare results from the Bayesian approach with other structure-preserving learning works. Finally, in Section \ref{sec: conclusions}, we provide concluding remarks and future research directions. 

\section{Background}
\label{sec: background}
In Section~\ref{sec: nonseparable}, we introduce nonseparable Hamiltonian systems, followed by their structure-preserving time integration using explicit symplectic integrators in Section~\ref{sec: explicit}. We then review the probabilistic modeling framework in Section~\ref{sec: bayesiansystemid}. This provides the necessary background for the structure-preserving Bayesian learning method described in Section~\ref{sec: learning}. 
\subsection{Nonseparable Hamiltonian systems}
\label{sec: nonseparable}
The governing equations for finite-dimensional canonical Hamiltonian systems are
\begin{equation}\label{eq:derivative}
      {\Rb{  \dot{\q}=\frac{\partial H(\q,\p)}{\partial \p}, \quad \quad \dot{\p}=-\frac{\partial H(\q,\p)}{\partial \q},}}
\end{equation}
where $H$, the Hamiltonian, is a function of the canonical position $\q$ and momentum $\p$.
Many Hamiltonian systems of interest to engineers and scientists (e.g., \cite{serra2019control}-\cite{colliander2010transfer}) are not additively separable with respect to functions of the position $\q$ and momentum $\p$. Such Hamiltonian systems are said to be \textit{nonseparable}. Unlike the separable Hamiltonian systems with $H(\q,\p)=T(\p) + U(\q) $, the governing equations in~\eqref{eq:derivative} cannot be further simplified for nonseparable Hamiltonians.

These governing equations possess physically meaningful geometric properties that can be described in the form of symmetries, symplecticity, first integrals, and energy conservation. Preservation of these qualitative features in a numerical simulation is crucial for accurate long-time prediction of Hamiltonian dynamics. Using ideas from geometric mechanics, the field of geometric numerical integration has developed a variety of structure-preserving time integrators for Lagrangian/Hamiltonian systems, e.g., \cite{hairer2006geometric,sharma2020review}.

\begin{example}[Cherry Problem] \label{cherry_example}
Consider the following nonseparable Hamiltonian from \cite{cherry1928v} with
    \begin{multline}\label{eq:cherry}
        H(q_1,q_2,p_1,p_2)=\frac{1}{2}(q_1^2 + p_1^2) -(q_2^2 + p_2^2) \\ + \frac{1}{2}p_2(p_1^2 - q_1^2) -q_1q_2p_1. 
    \end{multline}
This four-dimensional dynamical system is a challenging example because it possesses a negative energy mode (NEM) that leads to an explosive nonlinear growth of perturbations for arbitrarily small disturbances. These NEMs occur in several important infinite-dimensional dynamical systems, e.g., gravitational instability of interpenetrating fluids \cite{casti1998negative} and magnetosonic waves in the solar atmosphere \cite{joarder1997manifestation}.
\end{example}
\vspace{-0.75cm}
\subsection{Tao's explicit symplectic integrator}
\label{sec: explicit}

While explicit symplectic integration has been extensively studied for separable Hamiltonian systems and specific subclasses of nonseparable Hamiltonian systems, explicit symplectic approximations of general nonseparable Hamiltonian systems $H(\q,\p)$ is an active research topic. Recently, explicit symplectic integrators for arbitrary nonseparable Hamiltonian systems were presented in \cite{tao2016explicit}. The derivation of these integrators is based on the idea of an extended phase space. For an arbitrary nonseparable Hamiltonian $H(\q,\p)$, we first introduce fictitious configuration $\tilde{\q}$ and fictitious momentum $\tilde{\p}$ corresponding to $\q$ and $\p$, respectively.
Next, we define an augmented Hamiltonian 
\begin{multline}\label{eq:aug}
    \bar{H}(\q,\p,\tilde{\q},\tilde{\p}):=\underbrace{H(\q,\tilde{\p})}_{H_a} + \underbrace{H(\tilde{\q},\p)}_{H_b} \\ + \omega\cdot \underbrace{ \left( \Vert \q-\tilde{\q} \Vert_2^2/2 + \Vert \p-\tilde{\p} \Vert_2^2/2 \right)}_{H_c},
\end{multline}
where $H_a:=H(\q,\tilde{\p})$ and ${\Ra{H_b:=H(\tilde{\q},\p)}}$ correspond to two copies of the original nonseparable Hamiltonian system with mixed-up positions and momenta; $H_c$ is an artificial restraint;  and $\omega$ is a constant that controls the binding of the two copies. The governing equations for the augmented Hamiltonian system~\eqref{eq:aug} are  
{\small
 \begin{align*}\label{eq:derivative_augmented} 
        \dot{\q}&=\nabla_{\p}H(\tilde{\q},\p) + \omega (\p-\tilde{\p}), \quad \dot{\p}=-\nabla_{\q}H(\q,\tilde{\p}) - \omega (\q-\tilde{\q}), \\
        \dot{\tilde{\q}}&=\nabla_{\tilde{\p}}H(\q,\tilde{\p})  + \omega (\tilde{\p}-\p), \quad \dot{\tilde{\p}}=-\nabla_{\tilde{\q}}H(\tilde{\q},\p) - \omega (\tilde{\q}-\q).
\end{align*}
}%
The introduction of fictional variables $\tilde{\q}$ and $\tilde{\p}$ along with a specific choice for the augmented Hamiltonian in~\eqref{eq:aug} decouples the position $\q$ and $\p$ in the extended phase space, i.e., $\dot{\q}$ and $\dot{\p}$ are independent of $\q$ and $\p$, respectively. Unlike the original nonseparable Hamiltonian $H(\q,\p)$, the augmented Hamiltonian $\bar{H}(\q,\p,\tilde{\q},,\tilde{\p})$ is amenable to explicit symplectic integration.
 
In this work, we use second-order explicit symplectic method $\psi^{\dt}$ based on Strang splitting
\begin{equation}
\psi^{\Delta t}:=\psi^{\Delta t/2}_{H_a} \circ \psi^{\Delta t/2}_{H_b}\circ \psi^{\Delta t}_{\omega H_c} \circ \psi^{\Delta t/2}_{H_b}\circ \psi^{\Delta t/2}_{H_a},
\label{eq:exp_symp}
\end{equation} 
where $\psi^{\Delta t}_{H_a}, \psi^{\Delta t}_{H_b},$ and $\psi^{\Delta t}_{\omega H_c}$ are the time-$\dt$ flow of $H_a,H_b,$ and $\omega H_c$.{\Ra{ This allows us to obtain an explicit symplectic integrator for the augmented Hamiltonian $\bar{H}$ via composition of explicit symplectic Euler substeps with step size $\dt/2$.}} Explicit update equations for these individual flows can be written as
\begin{align*}
    \psi^{\Delta t}_{H_a}:& \begin{bmatrix} \q \\ \p \\\tilde{\q} \\\tilde{\p} \end{bmatrix} \to \begin{bmatrix} \q \\ \p -\Delta tH_{\q}(\q,\tilde{\p}) \\\tilde{\q} + \Delta t H_{\tilde{\p}}(\q,\tilde{\p}) \\\tilde{\p} \end{bmatrix}, \\ 
    \psi^{\Delta t}_{H_b}:& \begin{bmatrix} \q \\ \p \\\tilde{\q} \\\tilde{\p} \end{bmatrix} \to \begin{bmatrix} \q + \Delta t H_{\p}(\tilde{\q},\p)\\ \p  \\\tilde{\q}  \\\tilde{\p}-\Delta tH_{\tilde{\q}}(\tilde{\q},\p) \end{bmatrix}, \\
    \psi^{\Delta t}_{\omega H_c}:& \begin{bmatrix} \q \\ \p \\\tilde{\q} \\\tilde{\p} \end{bmatrix} \to \frac{1}{2}\begin{bmatrix}  \begin{pmatrix}
     \q + \tilde{\q} \\ \p + \tilde{\p}
     \end{pmatrix} + \mathbf R(\Delta t)  \begin{pmatrix}
     \q - \tilde{\q} \\ \p - \tilde{\p}
     \end{pmatrix} \\
      \begin{pmatrix}
          \q + \tilde{\q} \\ \p + \tilde{\p}
     \end{pmatrix} - \mathbf R(\Delta t)  \begin{pmatrix}
     \q - \tilde{\q} \\ \p - \tilde{\p}
     \end{pmatrix}
    \end{bmatrix},
\end{align*}
where $\mathbf R(\Delta t):=\begin{bmatrix} \cos (2\omega \Delta t) \mathbf{I} &  \sin (2\omega \Delta t) \mathbf{I} \\ -\sin (2\omega \Delta t) \mathbf{I} &  \cos (2\omega \Delta t) \mathbf{I} \end{bmatrix}$. 

Given a nonseparable Hamiltonian $H(\q,\p)$ with initial condition $[\q(0)^\top,\p(0)^\top]^\top=[\q_0^\top,\p_0^\top]^\top$, we obtain explicit symplectic approximations of the dynamics by integrating the augmented Hamiltonian $\bar{H}(\q,\p,\tilde{\q},\tilde{\p})$  with $[\q(0)^\top,\p(0)^\top,\tilde{\q}(0)^\top,\tilde{\p}(0)^\top]^\top=[\q_0^\top,\p_0^\top,\q_0^\top,\p_0^\top]^\top$.
\subsection{Bayesian system identification}
\label{sec: bayesiansystemid}
In the Bayesian learning framework, the Hamiltonians are unknown. Furthermore, they are parametrized in some approximation format, e.g., polynomial or neural network expansion.

We first model the system probabilistically as a hidden Markov model. Let $\X_k=[\q^\top_k, \p^\top_k ]^\top\in\reals^{2\dimx}$ be the random variable representing a belief of the system's state, and let the subscript $k$ be the index corresponding to the time $t_k$ at which output data $\bm{y}_k\in\reals^{\dimy}$ is observed. The output data is modeled as a realization of the random variable $\Y_k\in\reals^{\dimy}$ that represents a belief of the systems true output. 
The state-transition operator $\Psi:\reals^{2\dimx}\times\reals^{\dimpar}\rightarrow{\Ra{\reals^{2\dimx}}}$ maps the state $\X_k$ at an arbitrary time $t_k$ to the state $\X_{k+1}$, and the observation operator $h:\reals^{2\dimx}\times\reals^{\dimpar}\rightarrow\reals^{\dimy}$ maps the state $\X_k$ to the output $\Y_k$. Both operators may be parametrized by the random variable $\param\in\reals^{\dimpar}$. The model of the system dynamics is described by
\begin{equation}
\begin{aligned}\label{eq:system}
    \X_{k+1} &= \Psi(\X_k,\thetdyn) + \pmb{\xi}_k, &\quad \pmb{\xi}_k\sim\N(\bm{0},\Sigma(\thetsig)), \\
    \Y_{k} &= h(\X_k,\thetobs) + \pmb{\eta}_k, &\quad \pmb{\eta}_k\sim\N(\bm{0},\Gamma(\thetgam)).
\end{aligned}
\end{equation}
The parameters are partitioned as $\param=(\thetdyn,\thetobs,\thetsig,\thetgam)$; the process noise $\pmb{\xi}_k$ and measurement noise $\pmb{\eta}_k$ are additive, zero-mean Gaussian with unknown covariances $\Sigma(\thetsig)$ and $\Gamma(\thetgam)$, respectively. Gaussian noise is chosen not only because it is easy to work with, but also because for a given mean and covariance, the Gaussian distribution satisfies the Principle of Maximum Entropy~\cite{jaynes2003probability}, meaning it minimizes the amount of information that must be assumed beyond the first two moments.

In this model, we assume that there exists error/noise in both the state propagation and observation operators, and we therefore learn a stochastic system rather than a deterministic one. The source of the measurement error is typically sensor noise, and the term $\pmb{\eta}_k$ is included in nearly all system ID models. The error in the state propagation operator arises due to the uncertainty in both the parameters $\thetdyn$ and in the choice of the model structure. In other words, unless the model is perfectly known, the output of $\Psi$ will always carry some amount of error with respect to the ``true'' state. The process noise $\pmb{\xi}_k$ is not usually included in system ID models, but in~\cite{galioto2020bayesian}, it was shown to improve the robustness of system ID methods when faced with sparse and/or noisy data. In this paper, we will further demonstrate the benefits of accounting for the process error.

The goal of Bayesian inference is to compute a posterior distribution for $\param$ that represents available knowledge of the parameter values given all of the data that we have so far collected. This posterior distribution is denoted as $\probd(\param|\yn)$ where $\yn\coloneqq(\Y_1=\bm{y}_1,\ldots,\Y_{\numObs}=\bm{y}_{\numObs})$ is the collection of data points $\bm{y}_k$ treated as realizations of the random variables $\Y_k$ for $k=1,\ldots,\numObs$. Representing available knowledge as a distribution rather than a single estimated point is especially desirable in cases where the data are incomplete/weakly informative, such as sparse and/or noisy data. In these cases, the precision of the parameter estimate is limited, and it is therefore useful to be able to represent both the parameters and the predicted system output given these parameters as distributions.

To compute the posterior distribution, we factorize $\probd(\param | \yn)$ using Bayes' theorem:
\begin{equation}
    \probd(\param | \yn) = \frac{\like(\param;\yn)\probd(\param)}{\probd(\yn)}\propto\like(\param;\yn)\probd(\param),
\end{equation}
where $\like(\param;\yn)\coloneqq\probd(\yn | \param)$ is the likelihood, $\probd(\param)$ is the prior, and $\probd(\yn)$ is the evidence. Because $\probd(\yn)$ does not depend on the parameters, we treat it as a normalizing constant and do not need to evaluate it. The prior is chosen by the user, and is usually chosen to be easy to evaluate. The main challenge of evaluating the posterior then is the evaluation of the likelihood. Since the dynamics model is stochastic, the state $\X_k$ is also uncertain, and, although we are interested in the marginal likelihood $\like(\param;\yn)$, we are only able to directly compute the joint likelihood $\like(\xn,\param;\yn)$, where $\xn\coloneqq(\X_1,\ldots,\X_{\numObs})$. To extract the marginal likelihood from the joint likelihood requires the evaluation of the high-dimensional integral $\int\like(\xn,\param;\yn)\rm{d}\xn$. This costly computation can be avoided, however, by factorizing the marginal likelihood as follows: 
\begin{equation}\label{eq:like}
    \like(\param;\yn) = \prod_{k=1}^{\numObs}\like_k(\param;\yk),
\end{equation}
and by using the algorithm in~\cite[Thm 12.1]{sarkka2013bayesian} to recursively evaluate each term $\like_k(\param;\yn)$. When there exist nonlinearities in either of the functions $\Psi(\cdot,\cdot)$ or $h(\cdot,\cdot)$, the distributions $\like_k(\param;\yk)$ will be non-Gaussian, but we approximate them as Gaussian using the unscented Kalman filter (UKF) for computational efficiency. See~\cite{noh2019posterior, drovandi2022ensemble, khalil2015estimation} for a similar approach. {\Rb{The computational complexity of the integrator~\eqref{eq:exp_symp} is $\mathcal{O}(\dimx^3+\dimx G)$, where $G$ denotes the cost of evaluating the gradient of the Hamiltonian. Based on the analysis in~\cite{galioto2020bayesian}, the computational complexity of evaluating the marginal likelihood embedded with Tao's integrator is $\mathcal{O}(\numObs\dimx^3+\numObs\dimx G)$.}} For more details on this algorithm, see~\cite{galioto2020bayesian,galioto2020hamiltonian}. 

{\Rb{For other types of noise in system~\eqref{eq:system} such as multiplicative or non-Gaussian noise, different filtering algorithms can be used to evaluate the marginal likelihood with possibly improved approximation quality. If the noise is multiplicative, a Kalman filter variation designed for multiplicative noise~\cite{yang2002robust} can be used. If extreme non-Gaussianity is present in the estimation distribution, a particle filter~\cite{gustafsson2010particle} can be used to approximate the marginal likelihood with an empirical discrete distribution rather than a Gaussian, but the cost of the algorithm increases significantly. In Section~\ref{sec:results_multiple}, we break the assumption of additive, Gaussian noise and show that the Bayesian algorithm with the UKF can still deliver accurate estimates. The numerical results in Section~\ref{sec:results_multiple} demonstrate that the Bayesian algorithm has predictive power even when the assumptions of system~\eqref{eq:system} are not strictly met.}}

\section{Probabilistic Learning of Nonseparable Hamiltonian Systems}
\label{sec: learning}
This section presents a novel methodology for incorporating physical knowledge of nonseparable Hamiltonian systems into the Bayesian learning framework. This discussion includes the parametrization of the model dynamics in Section~\ref{sec: parameter}, the selection of a structure-preserving time integrator within the dynamics propagator in Section~\ref{sec: embed}, and an outline of the proposed computational procedure in Section~\ref{sec: data}.

\subsection{Parametrizing nonseparable Hamiltonian systems}
\label{sec: parameter}
Our goal is to leverage prior knowledge of the underlying system to inform the parametrization of a model. Recall from Section~\ref{sec: nonseparable} that a Hamiltonian system is defined by a function of the generalized position and momentum known as the Hamiltonian. To leverage this fact, we assume that the generalized coordinate frame is known and attempt to directly learn the Hamiltonian. By estimating a model of the Hamiltonian and deriving the dynamics from this model according to~\eqref{eq:derivative}, we guarantee that the predicted flow will also be Hamiltonian. This modeling choice makes the predictions more physically meaningful and also serves to restrict the model search space, which can help optimization.

As an example of a parametrization of the Hamiltonian, we consider a linear combination of nonlinear basis functions. 
In this work, we are concerned with nonseparable Hamiltonians, so we denote the model of the Hamiltonian generally as the differentiable function $\tilde{H}(\q,\p,\thetdyn)$.
Then we approximate the Hamiltonian up to some additive constant $C$ as
\begin{equation}\label{eq:ham_param}
    \tilde{H}(\q,\p,\thetdyn) = \Phi^\top(\q,\p)\thetdyn + C,
\end{equation}
where $\Phi(\q,\p)\in\reals^{\numBasis}$ is a vector whose $i$th component is the differentiable basis function $\phi_i:\reals^{2\dimx}\rightarrow\reals$ for $i=1,\ldots,\numBasis$. The gradient of the Hamiltonian follows
\begin{equation}
    \nabla \tilde{H}(\q,\p,\thetdyn) = \nabla\Phi^\top(\q,\p)\thetdyn,
\end{equation}
where $\nabla\Phi(\q,\p)\in\reals^{\numBasis\times2\dimx}$ is a matrix where the $i$th column is the gradient $\nabla\phi_i$ with respect to the state $[\q^\top, \p^\top]^\top$.

The learning framework of Section~\ref{sec: bayesiansystemid} is applicable to any arbitrary model parametrization including neural networks and other nonlinear approximations. 
\subsection{Embedding explicit symplectic integrators into the Bayesian learning framework }
\label{sec: embed}
In this subsection, we incorporate the stochastic nonseparable Hamiltonian into the Bayesian system ID framework. The key idea is to exploit knowledge of nonseparable Hamiltonian systems and their structure-preserving time integrators to inform the design of the  state propagator $\Psi(\q_k,\p_k;\param_{\Psi})$ in \eqref{eq:system}. We apply the second-order explicit symplectic integrator $\psi^{\Delta t}$ from \eqref{eq:exp_symp} to the parametrized nonseparable Hamiltonian \eqref{eq:ham_param} to obtain 
\begin{equation*}
    \psi^{\Delta t}_{\param}:=\psi^{\Delta t/2}_{\tilde H_a(\param)} \circ \psi^{\Delta t/2}_{\tilde H_b(\param)}\circ \psi^{\Delta t}_{\omega \tilde H_c} \circ \psi^{\Delta t/2}_{\tilde H_b(\param)}\circ \psi^{\Delta t/2}_{\tilde H_a(\param)},
\end{equation*}
where $\tilde H_a(\param):=\tilde H(\q,\tilde{\p},\param)$ and $\tilde H_b(\param):=\tilde H(\tilde{\q},\p,\param)$ correspond to two copies of the parametrized nonseparable Hamiltonian model with mixed-up positions and momenta, and $\omega \tilde H_c$ is an artificial restraint that controls the binding of the two parametrized copies. We use $\psi^{\Delta t}_{\param}$ to encode Hamiltonian structure into the propagator, i.e.,
\begin{equation*}
    \Psi(\q_k,\p_k;\param_{\Psi}):= L^{\dagger} \circ  \psi^{\Delta t}_{\param} \circ L,
\end{equation*}
where $L:[\q_k^\top, \p_k^\top]^\top \to [\q_k^\top, \p_k^\top, \q_k^\top, \p_k^\top]^\top$ duplicates the state to obtain the augmented state, and $L^{\dagger}:[\q_{k}^\top, \p_{k}^\top, \tilde{\q}_{k}^\top, \tilde{\p}_{k}^\top  ]^\top \to [\q_{k}^\top, \p_{k}^\top]^\top$ operates on the augmented state to yield the state in the original phase space.

This physics-informed choice of the propagator ensures that the learned model will  respect the geometric properties of the true process, i.e., the learned model will be a canonical nonseparable Hamiltonian system.
\subsection{Learning setups}
\label{sec: data}
We consider two primary learning setups: (1) data are collected from a single trajectory, and the goal is to predict the future state of this trajectory; and (2) data are collected from multiple trajectories with independent initial conditions, and the goal is to predict the state trajectory at an arbitrary initial condition. The first of these follows straightforwardly from the framework of Section~\ref{sec: bayesiansystemid}, but the second requires a small adaptation.

Let $\ynm$ be the collection of $\numObs$ data points measured from the trajectory with initial condition $\bm{x}_0^{(m)}$. The factorization of the likelihood~\eqref{eq:like} now becomes $\like(\param;\yn^{(1,\ldots,\numTraj)})=\prod_{m=1}^{\numTraj}\like(\param;\ynm)$, where each term $\like(\param;\ynm)$ is evaluated using~\eqref{eq:like}. This likelihood factorization has the added benefit of being easily parallelized for computational efficiency.

The second learning setup is also considered in other works related to learning Hamiltonian systems such as~\cite{wu2020structure, xiong2020nonseparable}. In this work, we focus specifically on comparing to~\cite{wu2020structure} since it uses the same parametrization~\eqref{eq:ham_param} that we consider. The method in~\cite{wu2020structure} attempts to learn a linear mapping from a dictionary of basis functions to the time derivatives by solving the following optimization problem:
\begin{equation}
    \min_{\vc\in\reals^{\numBasis}} \lVert \mA\vc - \vb \rVert_2,
\end{equation}
where each element of $\mA$ and $\vb$ is defined as
\begin{align}
    a_{ij} &= \frac{1}{K}\sum_{k=1}^{K}\Big(\nabla\phi_i(\q_k,\p_k)\cdot\nabla\phi_j(\q_k,\p_k)\Big), \\
    b_i & = \frac{1}{K}\sum_{k=1}^{K}\Big(\begin{bmatrix}\dot{\p}_k^\top&-\dot{\q}_k^\top\end{bmatrix}^\top\cdot\nabla\phi_i(\q_k, \p_k)\Big),
\end{align}
for $1\leq i, j\leq\numBasis$. Typically, data of $\dot{\q}$ and $\dot{\p}$ are not available and must therefore be numerically approximated. Following the approach of~\cite{wu2020structure}, we use a second-order finite difference method for this approximation. Since this method solves a linear least-squares problem, we henceforth refer to it as the least-squares (LS) method.

\subsection{Computational procedure}
\label{sec: computational}
This subsection details the computational procedure that we use to draw samples from the posterior and to make predictions using these samples.

To draw samples, we use Markov chain Monte Carlo (MCMC) sampling. Specifically, we use the delayed rejection adaptive Metropolis (DRAM)~\cite{haario2006dram} MCMC sampler. The proposal is a Gaussian random walk, and the covariance of this Gaussian is scaled by a factor of $\gamma=0.01$ to use as the second-tier proposal. We begin adaptation of the proposal covariance after 200 samples and add $\varepsilon=10^{-10}$ to the diagonal to maintain positive semidefiniteness. In all of the examples, we start sampling at the MAP point found using MATLAB's \texttt{fmincon} optimizer. 

Once we have samples from the posterior, we need to choose how to use them for prediction. Two specific values of the posterior predictive distribution, defined as $\probd(\mathcal{X}_{\numSteps}|\yn)=\int\probd(\mathcal{X}_{\numSteps}|\param)\probd(\param|\yn){\rm d}\param$ for arbitrary positive integer $T$, are used in the results. First, the maximum a posteriori (MAP) value $\param^{\text{MAP}}$ of the parameter posterior is simply its maximum. We use this to obtain the prediction $\bm{x}_k=\Psi^k(\bm{x}_0,\param^{\text{MAP}})$, where $\Psi^k$ denotes $k$ compositions of the state propagator $\Psi$ for any $k\in\mathbb{Z}_{+}$. 
The other value is the expected value of the posterior predictive distribution $\bm{x}_{1:\numSteps}=\mathbb{E}[\probd(\mathcal{X}_{\numSteps}|\yn)]$. The expected value is estimated as $\bm{x}_k=\frac{1}{\numMCMC}\sum_{i=1}^{\numMCMC}\Psi^k(\bm{x}_0,\param_i)$ for any $k\in\mathbb{Z}_{+}$, where $\param_i$ is the $i$th MCMC sample and $\numMCMC$ the total number of MCMC samples. Typically, we use a subset of samples drawn at regular intervals to perform this estimation for the sake of computational practicality. We will sometimes refer to this prediction point simply as the mean estimate. The decision of which point to use can depend on a number of variables such as the dynamics or the shapes of the parameter posterior and/or posterior predictive distributions. Since the points are easy to compute after sampling, heuristic choices are typically an acceptable decision method.
\section{Numerical Experiments}
\label{sec: numerical}
In this section, we apply the Bayesian learning method presented in Section \ref{sec: learning} to the Cherry problem discussed in Example~\ref{cherry_example}. First, we compare the structure-preserving Bayesian learning approach to a `physics-ignorant' approach where the propagator is equipped with a non-symplectic integrator of the same order that does not preserve the underlying Hamiltonian structure. The results of this comparison are provided in Section \ref{sec:results_single} and show the improvements in accuracy gained by using a structure-preserving symplectic integrator over a non-symplectic integrator of equivalent order accuracy. Section~\ref{sec:results_multiple}, gives a comparison between the proposed method and the LS method described in~\ref{sec: data}. This section demonstrates the greater robustness of the proposed approach to noisy data. 

In both sections, we parametrize the Hamiltonian with polynomials up to total order three for a total of $\numBasis = 34$ basis functions. The type of polynomials used in each example is given in the corresponding section. Each covariance matrix is parametrized as an identity matrix scaled by an unknown parameter: $\Sigma(\thetsig)=\theta_{35}\mathbf{I}_{4}$ and $\Gamma(\thetgam)=\theta_{36}\mathbf{I}_{4}$. We place a Laplacian prior on $\thetdyn$ to promote sparsity and half-normal priors on $\thetsig$ and $\thetgam$, which are the distributions of the absolute value of a zero-mean Gaussian random variable. The observation operator $h$ is the identity such that we measure the full state. The data and predictions are both generated using the explicit symplectic integrator using a time step of $\dt=0.01$.
\subsection{Training with a single initial condition}
\label{sec:results_single} 

In this study, we consider the case where noisy data from a single initial condition are used to learn the model, and we compare numerical results far outside the training data regime to  demonstrate the generalizability of the presented Bayesian learning algorithm.

For training the model, we select the initial condition $[q_1(0),q_2(0),p_1(0),p_2(0)]^\top=[0.15,0.1,-0.05,0.1]^\top$ and integrate the nonseparable Hamiltonian system up to $t=8$. We collect the data on the full state every $0.4$ time units for a total of 21 data points with zero-mean Gaussian noise with standard deviation of $\sigma=0.01$. For the learning problem, we parametrize the nonseparable Hamiltonian using monomial basis functions. To ensure proper convergence, we draw $2 \times 10^5$ samples, and discard the first $10^5$ samples as a conservative burn-in period.

\begin{figure*}[h]
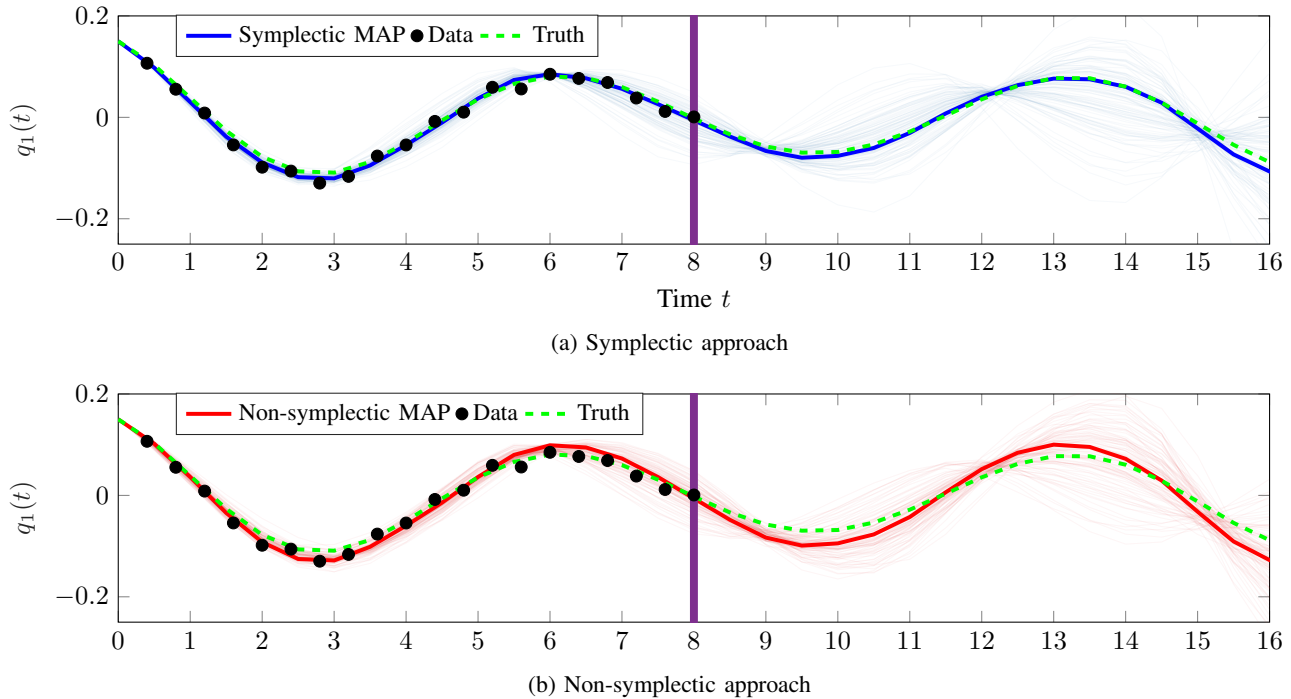

\centering
\begin{subfigure}{\textwidth}
          \setlength\fheight{6 cm}
          \setlength\fwidth{\textwidth}
\input{figures/cdc_exp.tex}
\caption{Symplectic approach}
\label{fig:single_exp}
    \end{subfigure} \\
\begin{subfigure}{\textwidth}
      \setlength\fheight{6 cm}
        \setlength\fwidth{\textwidth}
\input{figures/cdc_rk.tex}
\caption{Non-symplectic approach}
\label{fig:single_rk2}
    \end{subfigure} 
\caption{Reconstruction error comparison for the training with single IC case. The MAP estimate for the symplectic approach captures the nonseparable Hamiltonian dynamics accurately at $t=16$ which is $100 \%$ outside the training data whereas the MAP estimate for the non-symplectic approach yields inaccurate predictions outside the training data.  The posterior for both approaches grows as the system evolves further in time. The purple line indicates the end of the data collection period.}
 \label{fig:single}
\end{figure*}

Both symplectic and non-symplectic approaches in this study use a learning time step of $\Delta t=0.05$ (five times the time step used for prediction) to reduce computational time. Despite the large learning step, the symplectic approach gives accurate predictions far outside the training data regime, see Figure~\ref{fig:single_exp}. For the non-symplectic approach, we use a second-order Runge-Kutta integrator that does not preserve the underlying Hamiltonian structure.   

\begin{figure}[h]
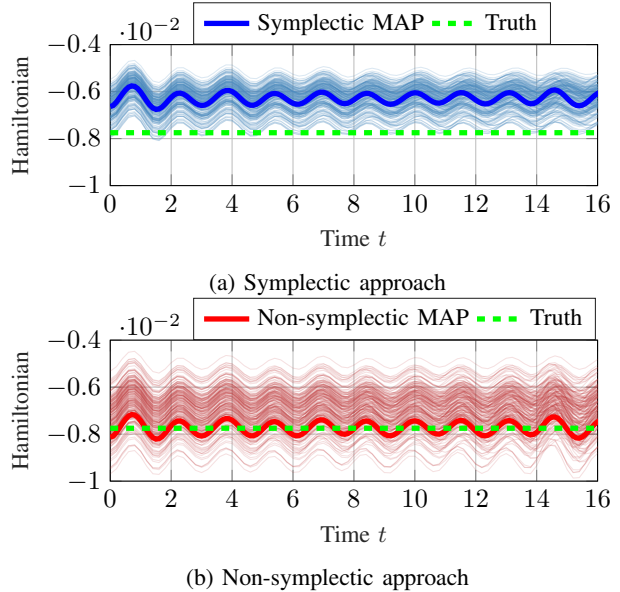

    \vspace{0.5cm}
    \begin{subfigure}{\linewidth}
           \setlength\fheight{3.5 cm}
           \setlength\fwidth{\textwidth}
           \input{figures/Hamiltonian_Exp.tex}
        \caption{Symplectic approach}
        \label{fig:H_Exp}
    \end{subfigure} \\
    \begin{subfigure}{\linewidth}
         \setlength\fheight{3.5 cm}
           \setlength\fwidth{\textwidth}
           \input{figures/Hamiltonian_RK.tex}
        \caption{Non-symplectic approach}
        \label{fig:H_RK}
    \end{subfigure}
    \caption{The learned Hamiltonians. The lighter blue and red lines in (a) and (b) are the samples from the posteriors from the symplectic and non-symplectic learning approaches. Both approaches estimate the true Hamiltonian accurately with an absolute energy error of less than $10^{-3}$. However, the non-symplectic approach has much greater uncertainty.}
    \label{fig:Hamcomparison}
\end{figure}

The MAP parameters of each posterior are used to simulate for $t=16$ (100 \% outside the training data regime) using the explicit symplectic integrator. 
Figure~\ref{fig:Hamcomparison} shows the learned Hamiltonian values for the MAP point and their corresponding posteriors for both symplectic and non-symplectic approaches. The true Hamiltonian is $-0.78 \times 10^{-2}$, and the value of the Hamiltonian learned by the symplectic approach is approximately $-0.68 \times 10^{-2}$. The Hamiltonian learned by the non-symplectic approach is approximately $-0.81 \times 10^{-2}$. Compared to the posterior from the symplectic approach  in Figure~\ref{fig:H_Exp}, the posterior from the non-symplectic approach in Figure~\ref{fig:H_RK} reflects greater uncertainty in the Hamiltonian estimate.

Figure~\ref{fig:single} compares the reconstructed trajectories of the system's first state using the MAP point and the samples from the posterior. For a fair comparison between the symplectic and non-symplectic approaches, each trajectory is integrated using the explicit symplectic integrator.\footnote{We assume that a non-symplectic learning approach would still actually integrate the learned model with a symplectic integrator. This setting is similar to the approach considered in \cite{greydanus2019hamiltonian} or \cite{zhong2019symplectic} where the authors use non-symplectic integrators for training structure-preserving neural networks and then integrate the learned Hamiltonian models using symplectic integrators.} We see in Figure~\ref{fig:single} that the MAP
estimate from the symplectic approach aligns closely with the truth for the entire length of the simulation, while the MAP point from the non-symplectic approach gives inaccurate results after the end of the training time interval at $t=8$. 

We also compute the following relative state error for both the training data regime and testing data regime
\begin{equation}\label{eq:rel_error}
    e(i:j) = \frac{\lVert\hat{\bm{x}}_{i:j}-\bm{x}_{i:j}\rVert_{F}}{\lVert\bm{x}_{i:j}\rVert_{F}}
\end{equation}
for positive integers $i\leq j$ where $\bm{x}_{i:j}\in\reals^{2\dimx\times (j-i+1)}$ is the true state trajectory at times $t_i,\ldots,t_j$, and $\hat{\bm{x}}_{i:j}$ is the estimate of this trajectory. 

In the training data regime ($t=0$ to $t=8$), the relative state error, $e(1:\numObs)$, for the symplectic approach is $0.0877$ whereas the relative state error for the non-symplectic approach is $0.1307$. In the testing data regime  ($t=8$ to $t=16$), the relative state error, $e(\numObs+1: 2\numObs)$, for the symplectic approach is $0.2173$ whereas the relative state error for the non-symplectic approach is $0.3326$. Thus, the symplectic approach achieves approximately $30 \%$ reduction in state error over the non-symplectic approach.
\subsection{Training with multiple initial conditions}
\label{sec:results_multiple}
In this subsection, we consider the case where data are collected from a number of different trajectories with independent initial conditions (ICs), and our goal is to estimate a model of the system that can give good prediction for an arbitrary IC. 

To gather the training data, we first randomly sample $\numTraj=5$ ICs from a Gaussian centered at the testing IC with standard deviation $0.05$. The testing IC is the same IC used in the previous example.  For each training IC, we integrate the system for up to $t=8$ and measure the state every $0.4$ time units for a total of 21 data points per trajectory including the IC. Any trajectories that diverged during this time period had their IC resampled. To complete the training dataset, we add 10\% relative noise, i.e., $\Y_k=\X_k(1+u_k)$ where $u_k\sim\mathcal{U}[-0.10, 0.10]$ and $\mathcal{U}$ is the uniform distribution. This is the same noise form used by~\cite{wu2020structure}.

We parametrize the Hamiltonian for the learning problem with Legendre polynomials as basis functions and integrate with the symplectic integrator with $\dt=0.01$ during training. To find the MAP, we use the LS estimate as the optimization starting point for $\thetdyn$ and the starting point $10^{-6}$ for variance parameters $\thetsig$ and $\thetdyn$. Then we draw $10^{4}$ samples from the posterior. 

\begin{figure*}
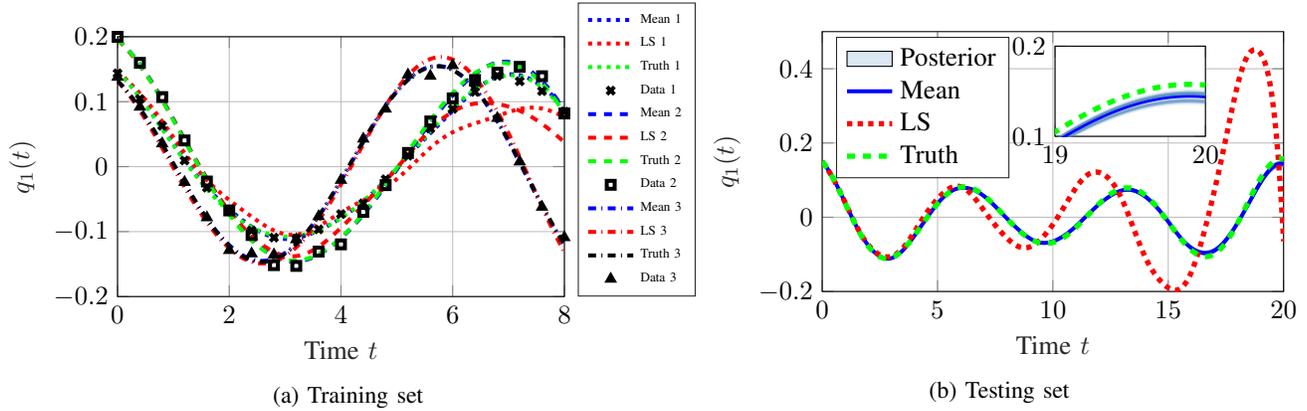

\vspace{0.5cm}
    \centering
    \begin{subfigure}{0.52\linewidth}
         \setlength\fheight{6.45 cm}
         \setlength\fwidth{\textwidth}
         \input{figures/q1_train.tex}
                \caption{Training set}
        \label{fig:q1_trainIC}
    \end{subfigure}
    \begin{subfigure}{0.44\linewidth}
         \setlength\fheight{6.45 cm}
         \setlength\fwidth{\textwidth}
         \input{figures/q1_test.tex}
        \caption{Testing set}
        \label{fig:q1_testIC}
    \end{subfigure}
    \caption{Left: The training data from 3 different ICs are shown alongside the point estimates from the LS~\cite{wu2020structure} and Bayesian methods. The Bayesian method is more robust to the noisiness/sparsity of the data. Right: The posterior predictive distribution and its mean are compared to the LS point on a trajectory outside the training set. The Bayesian method gives a good prediction over $t=20$, while the LS estimate deteriorates rapidly after about $t=5$.}
    \label{fig:LScomparison}
\end{figure*}

We compare the proposed structure-preserving learning method to the LS method and present the results in Figure~\ref{fig:LScomparison}, where every learned model was integrated with the explicit symplectic integrator. Note that although~\cite{wu2020structure} suggested a method for denoising the data, we found that using the raw data gave a better LS estimate in this problem. The training data and ground truth are shown alongside the point estimates from the LS and Bayesian methods in Figure~\ref{fig:q1_trainIC}. Even though the training ICs are all near the testing IC, we observe that the testing trajectory in Figure~\ref{fig:q1_testIC} differs significantly from the training trajectories. Figure~\ref{fig:q1_trainIC} also shows that the mean estimate is able to fit the training data fairly well, whereas the LS estimate struggles due to the noisiness/sparsity of the data. The fact that the mean estimate can fit the training data even when the measurement noise is non-Gaussian further demonstrates the robustness of the proposed approach. 

Next, we illustrate the posterior predictive distribution on the test IC in Figure~\ref{fig:q1_testIC} by plotting the trajectory estimate from every 200th MCMC sample for a total of 50 samples. We see that as time increases, the posterior starts to widen, reflecting the growing uncertainty in the state with time. We also observe that although both the LS and posterior predictive mean match the truth very closely through $t=4$, the LS estimate deviates from the truth much earlier and much more drastically than the mean estimate. Figure~\ref{fig:q1_testIC} demonstrates that the Bayesian method is able to find models that generalize well outside of the training set of the data.

To compare the methods quantitatively, consider the relative error~\eqref{eq:rel_error} over the first $k$ time steps, i.e., $e(1 : k)$. The relative error of the mean estimate on the test IC stays under 10\% for $t=18.22$, while the LS estimate stays under this threshold for only $t=1.49$. The Bayesian method therefore yields better long-time prediction performance.

\section{Conclusions}
\label{sec: conclusions}
We have developed a structure-preserving Bayesian learning approach for system ID of nonseparable Hamiltonian systems using stochastic dynamic models. This approach expands the Bayesian system ID work in~\cite{galioto2020bayesian} by incorporating knowledge of nonseparable Hamiltonian systems into the learning framework. The numerical results demonstrate the advantage of preserving the underlying geometric structure for a challenging nonseparable Hamiltonian system that possesses a negative energy mode. Hamiltonian models learned using the proposed method provide accurate predictions far outside the training time interval and perform well even for unknown initial conditions.

In future work, we would like to use neural networks for parametrization in the proposed framework for learning non-polynomial Hamiltonian systems. {\Rb{Moreover, while in this work we assume additive noises in the probabilistic model of the system dynamics~\eqref{eq:system}, in many applications the noise is multiplied by the state. It would be desirable to extend the current approach to dynamical systems with multiplicative noises. Finally, we have applied the structure-preserving Bayesian learning method to a system with four states. We would like to combine the proposed approach with Hamiltonian operator inference \cite{sharma2022hamiltonian} for structure-preserving learning of reduced-order models of large-scale Hamiltonian models.}}

\section{ACKNOWLEDGMENTS}
N. Galioto and A.A. Gorodetsky were funded by the AFOSR Computational Mathematics Program (P.M Fariba Fahroo). B. Kramer and H. Sharma were in part financially supported by the Ministry of Trade, Industry and Energy (MOTIE) and the Korea Institute for Advancement of Technology (KIAT) through the International Cooperative R\&D program (No.~P0019804, Digital twin based intelligent unmanned facility inspection solutions) and the Applied and Computational Analysis Program of the Office of Naval Research under award N000142212624.

\bibliography{references}
\bibliographystyle{IEEEtran}
\end{document}